\documentclass[12pt,notitlepage,twoside]{article}
\usepackage{amsmath,amssymb,amsthm}
\usepackage{graphpap,latexsym,epsf}
\usepackage{epsf}
\usepackage{color,graphics}

\usepackage{latexsym}
\usepackage{amsmath}
\usepackage{amsfonts}
\usepackage{amssymb}
\newcommand{\ristretto}{{\lower4pt\hbox{\bigg|}}}

\newcommand{\DHH}{\Delta_{\He^n} }
\newcommand{\ttt}{\widetilde}
\newcommand{\D }{\Delta }
\newcommand{\e }{\varepsilon }
\renewcommand{\l }{\lambda }
\newcommand\loc {{\text{\rm loc}}}
\newcommand{\n }{\nabla }
\newcommand{\nh }{\nabla_{\He^n} }

\newcommand{\ov}{\overline}
\newcommand\de {{\partial}}
\newcommand{\wtilde }{\widetilde}
\newcommand{\wk}{\rightharpoonup}
\newenvironment{pf}{\noindent{\bf Proof.}\enspace}{
\hfill$\Box$\medskip}
\newenvironment{pfn}[1]{\noindent{\bf Proof of {#1}.\enspace}}{
\hfill$\Box$\medskip}
\newcommand{\R}{\mathbb{R}}
\newcommand{\C}{\mathbb{C}}

\newcommand{\He}{\mathbb{H}}
\newcommand{\SC}{S^1_{\rm cyl}(\He^n)}

\newcommand{\conm}{{\mathop{\wk}\limits^{{\mathcal M}}}}
\newcommand{\SSS }{\mathbb{S}^{2n+1}}
\newcommand{\z }{\zeta }

\newtheorem{thm}{Theorem}[section]
\newtheorem{pro}[thm]{Proposition}
\newtheorem{lem}[thm]{Lemma}
\newtheorem{rem}[thm]{Remark}

\numberwithin{equation}{section}

\author{{\sc  Veronica Felli\thanks{S.I.S.S.A., Via Beirut 2-4, 34014
      Trieste, Italy; e-mail: \texttt{felli@sissa.it}}
       and Francesco Uguzzoni
\thanks{Dipartimento di Matematica, Universit\`{a} degli Studi di Bologna,
          Piazza di Porta S. Donato 5, 40126 Bologna,
          Italy;
          e-mail:  \texttt{uguzzoni@dm.unibo.it}
       }}}
\title { \Large \textbf{Some existence results \\ for
the Webster scalar curvature problem \\
in presence of symmetry}}

\begin{document}

\date{ }

\maketitle

{\footnotesize
\begin{abstract}

\noindent We prove some existence results for the
Webster scalar curvature problem on the Heisenberg
group and on the unit sphere of $\C^{n+1}$, under the
assumption of some natural symmetries of the
prescribed curvatures. We use variational and
perturbation techniques.



\end{abstract}

}

\section{Introduction}
In this paper we prove some existence results for the
equation
\begin{equation}\label{WH}
-\DHH u(\xi)=K(\xi)u(\xi)^{\frac{Q+2}{Q-2}},
    \qquad \xi\in \He^n,
\end{equation}
where $\DHH$ is the  sublaplacian on the Heisenberg
group $\He^n$ and $Q=2n+2$ is the homogeneous
dimension of $\He^n$. Our results provide existence of
solutions for the Webster scalar curvature problem on
 $\He^n$ and on the unit sphere
$\SSS$ of $\C^{n+1}$, under suitable assumption on the
prescribed
    curvatures. This problem is
the CR counterpart of the classical Nirenberg problem.
    In this paper we shall mainly assume that the prescribed
    curvature $K$ has a natural symmetry, namely a
    cylindrical-type symmetry.  Our main results are
    contained in Theorems \ref{t:Bianchi-Egnell-type},
    \ref{t:pert2}, \ref{t:pert3} and \ref{t:pert4} below.

 We remark that one of the main features of the above
equation \eqref{WH} is a lack of compactness due  both
to the criticality of the exponent $(Q+2)/(Q-2)$ and
to the unboundedness of the domain. Non-existence
results for \eqref{WH} can be obtained using the
Pohozaev-type identities of \cite{[GL]} under certain
conditions on $K$. In particular it turns out that a
positive solution $u$ to \eqref{WH} in the Sobolev
space $S^1_0(\He^n)$ (with the notation of Section 2)
satisfies the following identity
    $$\int_{\He^n}\langle (z,2t),\nabla K(z,t)\rangle
u(z,t)^{\frac{2Q}{Q-2}}dzdt=0$$ provided the integral
is convergent and $K$ is bounded and smooth enough.
This implies that there are no such solutions if
$\langle (z,2t),\nabla K(z,t)\rangle$ does not change
sign in $\He^n$ and $K$ is not constant.

The link between equation (\ref{WH}) and the Webster scalar curvature
problem on the sphere is briefly discussed below.
 Let us denote by $\theta_0$ the standard
contact form of the CR manifold $\SSS$. Given a smooth
function $\bar{K}$ on $\SSS$, the Webster scalar
curvature problem on $\SSS$ consists in finding a
contact form $\theta$ conformal to $\theta_0$ such
that the corresponding Webster scalar curvature is
$\bar{K}$ (for the definition of the Webster scalar
curvature see \cite{web}). This problem is equivalent
to solve the semilinear equation
\begin{equation}\label{WS}
b_n \Delta_{\theta_0} v(\z)+\bar{K}_{0}v(\z)
    =\bar{K}(\z)v(\z)^{b_n-1},
    \qquad \z\in \SSS,
\end{equation}
where $b_n=2+\frac{2}{n}$, $\Delta_{\theta_0}$ is the
sublaplacian on $(\SSS,\theta_0)$ and
$\bar{K}_0=\frac{n(n+1)}{2}$ is the Webster scalar
curvature of $(\SSS,\theta_0)$. If $v$ is a positive
solution to (\ref{WS}), then
$(\SSS,\theta=v^{\frac{2}{n}}\theta_0)$ has Webster
scalar curvature $\bar{K}$. Using the CR equivalence
$F$ (given by the Cayley transform, see definition
\eqref{F} below) between $\SSS$ minus a point and
$\He^n$, equation \eqref{WS} is equivalent to
\eqref{WH} with $K=\bar{K}\circ F^{-1}$, up to an
uninfluent constant.
     We refer to \cite{JL1} for a more
detailed presentation of the problem.

Indeed in the papers \cite{JL1, JL3, JL2}, Jerison and
Lee extensively studied the Yamabe problem on CR
manifolds (see also the recent papers \cite{Gamara,
GamYac}). On the contrary, at the authors' knowledge,
very few results have been established on the Webster
scalar curvature problem.
    In the recent paper \cite{mu} by Malchiodi
    and one of the authors, a
    new result is obtained in the
    perturbative case, i.e. when $K$ is assumed to be
    a small perturbation of a constant (see the
    papers \cite{agap, [ACPY3]},
    for analogous results concerning the
Riemannian context).

The aim of this paper is to begin to study a case
analogous to the radial one in the Riemannian setting.
The natural counterpart in our context seems to be
that of cylindrical curvatures
    $$K(z,t)=K(|z|,t)$$
     (see
Section 2 for all the notation) and not that of
``radial" ones $K=K(\rho)$. Indeed cylindrical
curvatures $K$ on $\He^n$ correspond on $\SSS$ to
curvatures $\bar{K}$ depending only on the last
complex variable of $\SSS\subseteq\C^{n+1}$,
    $$\bar{K}(\z_1,...,\z_{n+1})=\bar{K}(\z_{n+1}),$$  in
analogy with the Riemannian case where radial
curvatures $R$ on $\R^N$ correspond to curvatures
$\bar{R}(x_1,...,x_{N+1})=\bar{R}(x_{N+1})$ in
$\mathbb{S}^N\subseteq\R^{N+1}$, via the stereographic
projection.

However, the cylindrical case presents higher
difficulties with respect to the radial Riemannian
case. Indeed, when $K$ is cylindrical, one can reduce
equation \eqref{WH} to a two variables PDE, but not to
an ODE as in the radial case. Nevertheless we are able
to adapt a technique by Bianchi and Egnell \cite{be}
in order to obtain our first existence result Theorem
\ref{t:Bianchi-Egnell-type}. This technique consists
in a minimization on a space of cylindrically
symmetric functions and is based on a
concentration-compactness lemma which can be proved
just adapting the classical result by P. L. Lions
\cite{lions1,lions2} holding in the euclidean context.

We then deal with the perturbative case, obtaining
some results  via the abstract Ambrosetti-Badiale
finite dimensional reduction method \cite{ab1, ab2}
(see Theorems \ref{t:pert2}, \ref{t:pert3} and
\ref{t:pert4}). This method
allows to prove for the Webster scalar curvature the results found in
\cite{agap} for the scalar curvature problem in the
perturbative case. More precisely in \cite{agap}, Section 4, the
radial symmetry allows to reduce the perturbation
problem to the study of critical points of a one
variable function, thus obtaining more precise and
neat results than for the non radial case. Similarly
in our setting the cylindrical symmetry leads us to
treat a two variables problem and to find results like
Theorems \ref{t:pert2} and \ref{t:pert4}, which have
no counterpart in \cite{mu}, where the Webster
curvature problem is treated without requiring any
symmetry, and like Theorem \ref{t:pert3} which
requires assumptions of the type of \cite{mu} only on
the function $K$ restricted to the axis $\{z=0\}$. For
other  results related to the  Riemannian case in
presence of simmetry we refer to \cite{am,be2,be,hv}.

We finally remark that some other  results for
equation \eqref{WH} on the Heisenberg group have been
obtained in the  papers
 \cite{BRS, LW, U}. However,
     our
hypotheses on $K$ are very different from the ones in
such papers where $K$ is assumed to satisfy suitable
decaying conditions at infinity. In particular in
\cite{BRS} it is required an estimate of the type
$K_1(\rho) \DHH \rho \le K\le K_2(\rho) \DHH \rho$
($\rho$ is the homogeneous norm on $\He^n$ defined in
\eqref{rho} below) involving the degenerate term $\DHH
\rho$, which allows to ``radialize'' the problem and
to apply ODE methods.


\noindent {\it Acknowledgments.}  We wish to thank
Proff. A. Ambrosetti and E. Lanconelli for having
suggested the study of this problem. V. F. is
supported by M.U.R.S.T. under the national project
``Variational Methods and Nonlinear Differential
Equations''. F. U. is supported by University of
Bologna, funds for selected research topics.

\section{Notation and main results}
Denoting by $\xi=(z,t)=(x+iy,t)\equiv(x,y,t)$ the points of
$\He^n=\C^n\times\R\equiv \R^{2n+1}$, let us recall that the group law
on the Heisenberg group is
\[
(x,y,t)\circ(x',y',t')=(x+x',y+y',t+t'+2x'\cdot y-2x\cdot y')
\]
where $\cdot$ denotes the usual inner product in
$\R^n$. Let us denote by $\tau_{\xi}(\xi')=\xi\circ
\xi'$ the left translations, by $\delta_{\l}(\xi)=(\l
z,\l^2 t)$, $\l>0$ the natural dilations, by $Q=2n+2$
the homogeneous dimension, and by
\begin{equation}\label{rho}
\rho(\xi)=(|z|^{4}+t^{2})^{1/4}
\end{equation}
the homogeneous norm on $\He^n$. The Lie algebra of
left-invariant vector fields on $\He^n$ is generated
by
\[
X_j=\frac{\partial}{\partial x_j}+2 y_j\frac{\partial}{\partial
  t},\quad Y_j=\frac{\partial}{\partial y_j}-2x _j\frac{\partial}{\partial
  t},\quad j=1,2,\dots,n.
\]
The sub-elliptic gradient on $\He^n$ is given by $
\n_{\He^n}=(X_1,\dots,X_n,Y_1,\dots,Y_n) $ and the
Kohn Laplacian on $\He^n$ is the degenerate-elliptic
PDO
\[
\Delta_{\He^n}=\sum_{j=1}^n(X_j^2+Y_j^2).
\]

We will say that a function $f: \He^n\to\R$ is
continuous on $\overline{ \He}^n$ if it is continuous
and there exists $\lim_{(z,t)\to\infty}f(z,t)\in \R$.
In this case we
 will denote by $f(\infty)$ such a limit.

Let $K$ be a continuous function on
$\overline{\He}^n$. We shall always
suppose that $K$ has cylindrical symmetry, i.e.
$K(z,t)=\ttt{K}(|z|,t)$, and that $\ttt{K}$ is
  locally H\"older continuous in $]0,\infty[\times\R$.
 Let us consider the following
equation on $\He^n$
\begin{equation*}\tag{${\mathcal P}$}
-\D_{\He^n} u=K\, u^{Q^\star-1},\quad u>0 \quad\text{in}\ \He^n,
\end{equation*}
where $Q^{\star}=\frac{2Q}{Q-2}$. We will work in the space of
cylindrically symmetric functions of the Folland-Stein Sobolev space $S^1_0(\He^n)$, namely in
\[
\SC =\big\{u\in S^1_0(\He^n):\ u(z,t)=u(|z|,t)\big\},
\]
where $S^1_0(\He^n)$ is defined as the completion of
$C^{\infty}_0(\He^n)$ with respect to the norm
\[
\|u\|^2_{S^1_0(\He^n)}=\int_{\He^n}|\n _{\He^n}u|^2\, dz\,dt.
\]
Let us remark that $Q^\star$ is the critical exponent
for the embedding $S^1_0(\He^n)\hookrightarrow
L^{Q^\star}(\He^n)$. Choosing suitable regularization
functions, it is not difficult to recognize that we
have
\[ S^1_0(\He^n)=\bigg\{u\in
L^{Q^\star}(\He^n):\
\int_{\He^n}|\n_{\He^n}u|^2\,dz\,dt<\infty\bigg\}
\]
and that $\SC$ is equal to the closure in
$S_0^1(\He^n)$ of the set of cylindrically symmetric
$C_0^\infty$ functions. Let us also observe that $\SC$
is a Hilbert space endowed with the scalar product
$(u,v)=\int_{\He^n}\n _{\He^n}u\cdot \n _{\He^n}v\,
dz\,dt$. Let us denote by $S$ the best constant in the
Sobolev-type inequality (see \cite{JL3})
\begin{equation}\label{eq:sobolev}
S\|v\|^2_{Q^\star}\leq \|v\|^2_{S^1_0(\He^n)} \quad\forall\, v\in S_0^1(\He^n).\end{equation}
It is known (see \cite{JL3}) that all the positive  cylindrically
symmetric solutions in $S_0^1(\He^n)$ to the problem
\begin{equation*}
-\D_{\He^n}U= S^{\frac Q{Q-2}} U^{Q^\star -1}
\end{equation*}
are of the form
\begin{equation}\label{eq:ums}
U_{\mu,s}(z,t)=c_n \mu^{-\frac{Q-2}2} U_0\bigg(\frac r{\mu},\frac{t-s}{\mu^2}\bigg)
\end{equation}
(for $\mu>0$ and $s\in\R$) where  $r=|z|$,
\[
U_0(r,t)=\bigg(\frac{1}{t^2+(1+r^2)^2}\bigg)^{\!\!\frac{Q-2}4}
\]
and $c_n$ is a positive constant to be chosen in such
a way that $\int_{\He^n}|\n_{\He^n} U_{\mu,s}|^2=1$.
By solutions to problem $({\mathcal P})$ we mean weak
solutions in the sense of $S^1_0(\He^n)$. On the other
hand, under our hypotheses we have that solutions in
the $\SC$-sense are also solutions in the
$S_0^1(\He^n)$-sense, as it is shown in Lemma
\ref{wcs} in the Appendix.

Our first result is the counterpart in the Heisenberg context of a
result of Bianchi and Egnell \cite{be} about radial solutions of the
corresponding problem for the Laplacian on the euclidean space.
\begin{thm}\label{t:Bianchi-Egnell-type}
Assume that $K$ is a continuous cylindrically symmetric function on
$\He^n$ such that there exists
$K(\infty)=\lim_{(z,t)\to\infty}K(z,t)\in\R$, $K$
is positive somewhere and
\begin{equation}\label{eq:1}
\sup_{t\in\R}K(0,t)\leq K(\infty).
\end{equation}
If either $K(\infty)\leq0$ or there exist $\mu>0$ and $s\in\R$ such that
\begin{equation}\label{eq:2}
\int_{\He^n}(K(z,t)-K(\infty))U_{\mu,s}^{Q^\star}\,dz\,dt\geq0,
\end{equation}
then problem $({\mathcal P})$
has a cylindrically symmetric solution.
\end{thm}
\begin{rem}
The solution $u$ found in Theorem
\ref{t:Bianchi-Egnell-type} above, satisfies the decay
condition $u=O(\rho^{2-Q})$ at infinity (see e.g.
\cite[Proposition 2]{mu}). Moreover, by means of
standard regularization techniques based on the
results of Folland and Stein \cite{FS}, one can prove
that $u$ is smooth if $K$ is smooth. We finally remark
that Theorem \ref{t:Bianchi-Egnell-type} gives also an
existence result for the Webster scalar curvature
problem on the sphere $\SSS$, by means of the CR
equivalence $F:\SSS\setminus\{(0,...,0,-1)\}\to\He^n$,
\begin{equation}\label{F}
F(\z_1,...,\z_{n+1})
    =\bigg(\frac{\z_1}{1+\z_{n+1}}
    ,...,\frac{\z_n}{1+\z_{n+1}},{\rm {Re}}\Big(i\frac{1-\z_{n+1}}{1+\z_{n+1}}\Big)
    \bigg).
\end{equation}
We also remark that the set
$\{(0,t)|t\in\R\}\subseteq\He^n$ (i.e. the center of
the Heisenberg group) corresponds via $F$ to the
circle
$\{(0,w)\in\C^n\times\C|w\in{\mathbb{S}}^1\}\subseteq\SSS$.
Hence condition \eqref{eq:1} is equivalent to
$$\max_{w\in{\mathbb{S}}^1} {\bar{K}(0,w)} \leq\bar{K}(0,-1)$$
(where $\bar{K}=K\circ F$ is the prescribed curvature
on $\SSS$).
\end{rem}

Theorem \ref{t:Bianchi-Egnell-type} is proved by a minimization
technique which makes use of some concentration compactness argument;
since the solution is found as a constrained minimum, it is not
possible to obtain in the same way an analogous result with
 inequalities (\ref{eq:1}) and~(\ref{eq:2}) in the opposite sense.

In the second part of the paper, we shall deal with
the case in which $K$ is close to a constant, namely
$K(z,t)=1+\e k(z,t)$. We will consider the
perturbation problem
\begin{equation*}\tag{${\mathcal P}_{\e}$}
-\D_{\He^n} u=(1+\e k)\, u^{Q^\star-1},\quad u>0 \quad\text{in}\ \He^n,
\end{equation*}
where $\e$ is a small parameter and $k$ is a bounded
cylindrically symmetric function on~$\He^n$. Following the Ambrosetti and
Badiale \cite{ab1,ab2} finite dimensional reduction method we are able
to prove some perturbative existence results, which in most cases require
weaker assumptions than Theorem \ref{t:Bianchi-Egnell-type}. Our first
perturbation result is the counterpart of \cite{agap}, Theorem 4.5.
\begin{thm}\label{t:pert2}
Assume that $k$ is a continuous
cylindrically symmetric function on $\He^n$ with $k(\infty):=\lim_{(z,t)\to\infty}k(z,t)\in\R$
and that there exist $\mu>0$ and $s\in\R$ such that either
\begin{equation}\label{eq:5}
\int_{\He^n}\Big(k(z,t)-\sup_{\sigma\in\R}k(0,\sigma)\Big)U_{\mu,s}^{Q^\star}\,dz\,dt>0
\end{equation}
or
\begin{equation}\label{eq:oppositesense}
\int_{\He^n}\Big(k(z,t)-\inf_{\sigma\in\R}k(0,\sigma)\Big)U_{\mu,s}^{Q^\star}\,dz\,dt<0.
\end{equation}
Then problem $({\mathcal P}_{\e})$ has a cylindrically
symmetric solution for $|\e|$ sufficiently small.
\end{thm}
\begin{rem}\label{r:comparison}
{}\par\noindent\begin{enumerate}
\item In fact, from the proof it will be clear that we do not need to
  assume that the limit of $k$ at $\infty$ exists, if we assume instead that either
\begin{align*}
\int_{\He^n}\!\!\!\big(k(z,t)-\limsup_{(\mu,s)\to\infty}k(\mu,s)\big)U_{\mu,s}^{Q^\star}\,dz\,dt>0
\end{align*}
or
\begin{align*}
\int_{\He^n}\!\!\!\big(k(z,t)-\liminf_{(\mu,s)\to\infty}k(\mu,s)\big)U_{\mu,s}^{Q^\star}\,dz\,dt<0.
\end{align*}
\item Note the presence of the strict inequality in
(\ref{eq:5}) and (\ref{eq:oppositesense}) which is due to technical reasons.
\item  Assumption (\ref{eq:5}) is weaker
than assumptions~(\ref{eq:1}) and (\ref{eq:2}), with
strict inequality in at least one of them. The case covered by
(\ref{eq:oppositesense}) has no counterpart in Theorem
\ref{t:Bianchi-Egnell-type}. Actually in such a case the solutions we
will found are not constrained minima like the ones found in Theorem
\ref{t:Bianchi-Egnell-type}.
\end{enumerate}
\end{rem}
Using the perturbation method, it is also possible to find some other
results, requiring assumptions on the behavior of $k$ on the axis
$\{z=0\}$ instead of integral assumptions of the type (\ref{eq:5}) and
(\ref{eq:oppositesense}). The following result is the analogous for
the Heisenberg group of \cite{agap}, Theorem 4.4.
\begin{thm}\label{t:pert3}
Assume that $k$ is a cylindrically
symmetric function such that $\bar k=k\circ F$ is a smooth function on
${\mathbb S}^{2n+1}$, and that there exists a point
$(0,\bar s)\in\He^n$ such that
$k(0,\bar s)=\max_{\sigma}k(0,\sigma)$ and
\[
\Delta_{x,y}k(0,\bar s)>0
\]
where
$
\Delta_{x,y}
  k=\sum_{i=1}^n[\frac{\partial^2 k}{\partial
  x_i^2}+\frac{\partial^2 k}{\partial y_i^2}]
$.
 Then problem $({\mathcal P}_{\e})$ has a cylindrically symmetric solution for $|\e|$
sufficiently small.
\end{thm}
\begin{rem}\label{r:lapl}
It is also possible to find solutions under the assumption that there
exists a point $(0,\bar s)\in\He^n$ such that
$k(0,\bar s)=\min_{\sigma}k(0,\sigma)$ and
$\Delta_{x,y}k(0,\bar s)<0$.
\end{rem}
\begin{rem}\label{r:boun}
Let us remark that the assumption that $k$ comes from
a regular function on the sphere through the Cayley transform implies
that $k$ has finite limit at $\infty$ and that $|k|, |\n k|, |\n^2
k|$, and $|t^2\partial^2_t k|$ are bounded.
\end{rem}\noindent
Our last result is inspired by \cite{am}, Theorem 5.1(a).
\begin{thm}\label{t:pert4}
Assume that $k$ is a cylindrically symmetric continuous function such
that there exists $\lim_{(z,t)\to\infty}k(z,t)$ and  $k(0,t)\equiv k(\infty)$ for any
  $t\in\R$. Then problem $({\mathcal P}_{\e})$ has a cylindrically symmetric solution for $|\e|$ sufficiently small.
\end{thm}

\section{Proof of Theorem \ref{t:Bianchi-Egnell-type}}

In this section we shall prove Theorem
\ref{t:Bianchi-Egnell-type} by finding a solution of
$({\mathcal P})$ as a minimizer on the constraint
\[
M=\left\{ f\in \SC :\ \|f\|_{S^1_0(\He^n)}=1,\ \
\int_{\He^n}K\,|f|^{Q^\star}>0\right\}.
\]
Note that the assumption that $K$ is positive somewhere ensures that
$M$ is nonempty. Let us consider the minimum problem
\begin{equation*}\tag{${\mathcal I}_K$}
\gamma=\inf_{u\in M} {\cal F}_K(u)
\end{equation*}
where ${\cal F}_K: M\to\R$ is defined by
\[
 {\cal F}_K(u)=\left(\int_{\He^n}K\,|u|^{Q^\star}\right)^{-\frac Q{Q^\star}}=
\left(\int_{\He^n}K\,|u|^{Q^\star}\right)^{-\frac {Q-2}{2}}.
\]
Suppose that $u\in M$ attains the infimum in
$({\mathcal I}_K)$. Then $u$ is a critical point of
${\cal F}_K$ constrained on $M$. Then there exists a
Lagrange multiplier $\l\in\R$ such that
\[
( {\cal F}_K'(u),v)=\l ({\mathcal G}'(u),v)\quad\forall\, v\in \SC
\]
where ${\mathcal G}(u)=\int_{\He^n}|\n_{\He^n} u|^2$, namely for any $v\in \SC$
\begin{equation}\label{eq:3}
-Q\left(\int_{\He^n}K\,|u|^{Q^\star}\right)^{-\frac
 Q{2}}\int_{\He^n}K\,|u|^{Q^\star-2}u\,v=2\l
 \int_{\He^n}\n_{\He^n} u\cdot \n_{\He^n} v.
\end{equation}
Testing (\ref{eq:3}) with $v=u$, we can easily compute the value of
$\l$, thus finding $\l=-\frac Q2  {\cal F}_K(u)$. Hence (\ref{eq:3}) can be written in the form
\[
\int_{\He^n}\n_{\He^n} u\cdot \n_{\He^n}
v=\left(\int_{\He^n}K\,|u|^{Q^\star}\right)^{\!\!-1}\int_{\He^n}K\,|u|^{Q^\star-2}u\,v\quad\forall\,
v\in \SC
\]
which is equivalent to the fact that $u$ is a weak
solution in the $\SC$-sense of the equation
\[
-\D_{\He^n}u=\left(\int_{\He^n}K\,|u|^{Q^\star}\right)^{\!\!-1}K\,|u|^{Q^\star-2}u.
\]
Without loss of generality, we can assume $u\geq0$; otherwise one
takes $|u|$ after noticing that if $u\in M$, then $|u|\in M$ and $
{\cal F}_K(|u|)= {\cal F}_K(u)$.
    From Lemma \ref{wcs}, \cite[Proposition~2]{mu} and
    the Harnack inequality proved in \cite[Proposition 5.12]{JL1}
    it follows that $u$ is strictly positive and hence
$u$ satisfies
\[
-\D_{\He^n}u=\left(\int_{\He^n}K\,|u|^{Q^\star}\right)^{-1}K\,u^{Q^\star-1}.
\]
Let us consider the rescaled function
\[
\bar u=\left(\int_{\He^n}K\,|u|^{Q^\star}\right)^{\frac 1{2-Q^\star}}u.
\]
It is easy to check that $\bar u>0$ satisfies
\[
-\D_{\He^n}\bar u=K\,\bar u^{Q^\star-1},
\]
hence $\bar u$ is a solution to problem $({\mathcal
P})$. Therefore, the above argument shows that the
existence of minimizers of $({\mathcal
  I}_K)$
provides a weak solution in $\SC$ to $({\mathcal P})$
(and hence in~$S^1_0(\He^n)$, thanks to Lemma
\ref{wcs}).\par A sufficient condition for the
existence of minimizers of $({\mathcal
  I}_K)$ is given in the following lemma.
\begin{lem}\label{l:1}
Let $\delta=\sup\{K^+(0,t):\ t\in\R\}$ and assume that
there exists $V\in M$ such that
\[
 {\mathcal F}_K(V)\leq \inf_{v\in M}{\mathcal F}_{\delta}(v)={\mathcal
 F}_{\delta}(U_{\mu,s})=\delta^{-\frac{Q-2}2}S^{\frac Q2}.
\]
If $\delta=0$ it is enough to assume $M\not =\emptyset$ (i.e. $K$ positive
somewhere). Then the infimum in $({\mathcal I}_K)$ is attained.
\end{lem}
\begin{pf}
Let $(u_m)_m$ be a minimizing sequence for $({\mathcal I}_K)$, i.e.
\[
\left(\int_{\He^n}K\,|u_m|^{Q^\star}\right)^{\frac{2-Q}{2}}\mathop{\longrightarrow}\limits_{m\to+\infty}
\gamma,\quad \int_{\He^n}K\,|u_m|^{Q^\star}>0,\quad \int_{\He^n}|\n_{\He^n} u_m|^2=1.
\]
Up to a subsequence, we can assume that $u_m\wk u$ in
$\SC$ (and in $S^1_0(\He^n)$ and
$L^{Q^\star}(\He^n)$). In view of Lemma \ref{l:A1} we
have that $u_m$ converges to $u$ in $L^q(C)$ for any
set $C$ of the type $\{z: 0<c_1\leq |z|\leq
c_2\}\times [-c_3,c_3]$ (and hence on any compact set
away from the axis $\{z=0\}$) and for any $1\leq
q<+\infty$. From a suitable $\He^n$-version of the
concentration-compactness principle of P. L. Lions
(see Theorem \ref{t:cc}) we have that there exist some
nonnegative finite regular Borel measures $\nu,\nu_K,
\mu$ on $\overline{\He}^n$, an at most countable index
set $J$, a sequence $(z_j,t_j)\in\He^n$,
$\nu^j,\nu^{\infty}\in(0,\infty)$ such that (passing
to a subsequence) the following convergences in the
weak sense of measures hold
\begin{align}
|u_m|^{Q^\star}\,dx\   &\conm\ \nu=|u|^{Q^\star}+\sum_{j\in J}\nu^j\delta_{(z^j,t^j)}+\nu^{\infty}\delta_{\infty}\label{eq:wk1}\\
K|u_m|^{Q^\star}\,dx\   &\conm\ \nu_K=K|u|^{Q^\star}+\sum_{j\in J}K(z^j,t^j)\nu^j\delta_{(z^j,t^j)}+K(\infty)\nu^{\infty}\delta_{\infty}\label{eq:wk2}\\
|\n_{\He^n}u_m|^2\,dx\ &\conm\ \mu\geq |\n_{\He^n}u|^2+\sum_j S (\nu^j)^{\frac
  2{Q^\star}}\delta_{(z^j,t^j)}+S(\nu^{\infty})^{\frac
  2{Q^\star}}\delta_{\infty\label{eq:wk3}}.
\end{align}
Moreover if $(z^j,t^j)$ does not belong to the axis $\{z=0\}$, i.e. if
$z^j\not=0$, then $(z^j,t^j)$ is in some compact set of the type
$\{z:0<c_1\leq |z|\leq c_2\}\times [-c_3,c_3]$. Hence, in view of
Lemma \ref{l:A1} in the Appendix, $\int_C|u_m-u|^{Q^\star}\to 0$ and
consequently $\int_C|u_m|^{Q^\star}\to \int_C|u|^{Q^\star}$. Take some
continuous nonnegative function $\varphi$ with compact support and
satisfying $\varphi(z^j,t^j)\not=0$. For such a $\varphi$ we have that
$\int_C \varphi |u_m|^{Q^\star}\to \int_C\varphi |u|^{Q^\star}$ and hence
\[
0= \nu^j\varphi (z^j,t^j).
\]
Therefore it must be $\nu^j=0$, so that we can assume, without loss of
generality, that $z^j=0$ for any $j\in J$. We claim that $u\in M$. Let
us distinguish two cases.
\\[.5cm]
{\bf Case \boldmath $\delta>0$.} We can assume $\gamma<\delta^{-\frac{Q-2}2} S^{\frac
    Q2}$. Otherwise if $\gamma\geq \delta^{-\frac{Q-2}2} S^{\frac
    Q2}$ we have ${\cal F}_k(V)\leq \delta^{-\frac{Q-2}2} S^{\frac
    Q2}\leq \gamma$ and hence $V$ is a minimizer. Let $\gamma<\delta^{-\frac{Q-2}2} S^{\frac
    Q2}$. Since  $u_m$ is a minimizing sequence, from (\ref{eq:wk2})
  and (\ref{eq:wk3}) we have
\begin{align}
\delta S^{-\frac {Q^\star}2}
<\gamma^{-\frac
2{Q-2}}&=\nu_K
(\ov{\He^n})=\int_{\He^n}K|u|^{Q^\star}+\sum_j
K(0,t^j)\nu^j+K(\infty)\nu^{\infty},\label{sette}\\
1&\geq \int_{\He^n}|\n_{\He^n} u|^2+\sum_j S (\nu^j)^{\frac
  2{Q^\star}}+S(\nu^{\infty})^{\frac
  2{Q^\star}}\label{otto}.
\end{align}
We claim that $\int_{\He^n}|\n_{\He^n} u|^2=1$. By a way of
contradiction, assume that
$\int_{\He^n}|\n_{\He^n} u|^2=\rho$ with $\rho\in[0,1[$. From (\ref{otto}) we have
\begin{equation}\label{nove}
S^{-\frac{Q^\star}2}\geq (1-\rho)^{-\frac{Q^\star}2}\Big[
  \sum_j(\nu^j)^{\frac 2{Q^\star}}+(\nu^{\infty})^{\frac
  2{Q^\star}}\Big]^{\frac {Q^\star}2}.
\end{equation}
(\ref{sette}) and (\ref{nove}) imply that
\[
\gamma^{-\frac 2{Q-2}}>\delta (1-\rho)^{-\frac{Q^\star}2}\Big[
  \sum_j(\nu^j)^{\frac 2{Q^\star}}+(\nu^{\infty})^{\frac
  2{Q^\star}}\Big]^{\frac {Q^\star}2}
\]
hence $\gamma^{-\frac 2{Q-2}}>\delta (1-\rho)^{-\frac{Q^\star}2}\big(
  \sum_j\nu^j+\nu^{\infty}\big)$. Therefore
\begin{equation}\label{dieci}
\delta\Big(
  \sum_j\nu^j+\nu^{\infty}\Big)<(1-\rho)^{\frac{Q^\star}2} \gamma^{-\frac 2{Q-2}}.
\end{equation}
(\ref{sette}), (\ref{dieci}), and the definition of $\delta$ imply that
\begin{equation}\label{undici}
\gamma^{-\frac 2{Q-2}}=\int_{\He^n}K|u|^{Q^\star}+\sum_j
K(0,t^j)\nu^j+K(\infty)\nu^{\infty}<\int_{\He^n}K|u|^{Q^\star}+
(1-\rho)^{\frac{Q^\star}2} \gamma^{-\frac 2{Q-2}}.
\end{equation}
Hence $\int_{\He^n}K|u|^{Q^\star}>0$. Set $\bar u=
u\rho^{-1/2}\in M$; we have
\begin{equation}\label{dodici}
\frac{\int_{\He^n}K|u|^{Q^\star}}{\rho^{Q^\star/2}}=\int_{\He^n}K|\bar
u|^{Q^\star}<\gamma^{-\frac  2{Q-2}}.
\end{equation}
>From (\ref{undici}) and (\ref{dodici}) we obtain that
\[
\gamma^{-\frac  2{Q-2}}<\gamma^{-\frac
  2{Q-2}}\Big[\rho^{\frac{Q^\star}2}+(1-\rho)^{\frac{Q^\star}2}\Big]
\]
namely $1<\rho^{{Q^\star}/2}+(1-\rho)^{{Q^\star}/2}\leq
\rho+(1-\rho)=1$ which is not possible.
\\[.5cm]
{\bf Case \boldmath $\delta=0$.} In this case $K(\infty)\leq0$ and
$K(0,t)\leq 0$ for any $t\in\R$ and thus (\ref{sette}) and (\ref{dodici}) imply
\[
0<\gamma^{-\frac2{Q-2}}\leq \int_{\He^n}K|u|^{Q^\star}\leq
\gamma^{-\frac2{Q-2}}\bigg(\int_{\He^n}|\n_{\He^n}u|^2\bigg)^{\frac {Q^\star}2}.
\]
Hence $\int_{\He^n}|\n_{\He^n}u|^2\geq 1$. Since from
(\ref{otto}) we have $\int_{\He^n}|\n_{\He^n}u|^2\leq 1$, we can conclude $\int_{\He^n}|\n_{\He^n}u|^2= 1$.
\\[.5truecm]
The claim that $u\in M$ is thereby proved. From (\ref{otto}) we obtain
also that it must be $\nu^j=\nu^{\infty}=0$. Hence, in view of
(\ref{sette}), we deduce that the minimum in $({\mathcal I}_K)$ is
attained by $u$.
\end{pf}

\begin{pfn}{Theorem \ref{t:Bianchi-Egnell-type} completed} Let $\delta$ be as in Lemma \ref{l:1}. If
$K(\infty)\leq0$ then $\delta=0$. If $K(\infty)>0$ then
$\delta=K(\infty)>0$ and from (\ref{eq:2}) we have
\[
\int_{\He^n}K U_{\mu,s}^{Q^\star}\geq\delta\int_{\He^n} U_{\mu,s}^{Q^\star}
\]
and hence
\[
\bigg(\int_{\He^n}K U_{\mu,s}^{Q^\star}\bigg)^{\frac{2-Q}2}\leq
\bigg(\int_{\He^n}\delta U_{\mu,s}^{Q^\star}\bigg)^{\frac{2-Q}2}.
\]
Lemma \ref{l:1} (with $V=U_{\mu,s}$) allows us to
conclude.
\end{pfn}


\section{The perturbation problem}

In this section, we focus our attention on the case in which $K$ is
close to a constant, namely $K(z,t)=1+\e k(z,t)$. We deal with the perturbation
problem
\begin{equation*}\tag{${\mathcal P}_{\e}$}
-\D_{\He^n} u=(1+\e k)\, u^{Q^\star-1},\quad u>0 \quad\text{in}\ \He^n,
\end{equation*}
where $\e$ is a small real perturbation parameter and
$k$ is a bounded cylindrically symmetric function on
$\He^n$. Our approach is based on the finite
dimensional reduction method developed by Ambrosetti
and Badiale \cite{ab1,ab2} and recently applied by
Malchiodi and one of the authors \cite{mu} to the
problem of prescribing the Webster scalar curvature on
the unit sphere of $\C^{n+1}$. Cylindrically symmetric
solutions of $({\mathcal P}_{\e})$ can be obtained as
critical points on the space $\SC$ of the functional
\[
f_{\e}(u)=\frac12
\int_{\He^n}|\n_{\He^n}u|^2-\frac1{Q^{\star}}\int_{\He^n}(1+\e
k)u_+^{Q^\star}.
\]
Indeed, if $u$ is a nontrivial critical point of
$f_{\e}$, testing $f_{\e}'(u)$ with $u_-=\max\{-u,0\}$
we obtain that
$0=(f_{\e}'(u),u_-)=-\|u_-\|_{S^1_0(\He^n)}^2$  and
hence $u_-=0$.
    From Lemma \ref{wcs}, \cite[Proposition
2]{mu} and
    the Harnack inequality proved in \cite[Proposition 5.12]{JL1}
    it follows that~$u>0$.\par For $\e=0$,
the unperturbed functional $f_0$ has a manifold of
critical points $Z$ given by
\begin{equation}\label{eq:critman}
Z=\big\{z_{\mu,s}:\ \mu>0,\ s\in\R\big\}
\end{equation}
where $z_{\mu,s}$ is given, up to a constant, by the function
$U_{\mu,s}$ defined in (\ref{eq:ums}), namely
\begin{equation*}
z_{\mu,s}(z,t)= \mu^{-\frac{Q-2}2} \omega\bigg(\frac
{z}{\mu},\frac{t-s}{\mu^2}\bigg)
\end{equation*}
where
\begin{equation}\label{eq:defomega}
\omega(z,t)=(Q-2)^{\frac{Q-2}2}\big(t^2+(1+|z|^2)^2\big)^{-\frac{Q-2}4}.
\end{equation}

\subsection{The abstract perturbation method}\label{sec:abs}

For the reader's convenience, here we recall the
abstract result we will use in the sequel, for the
proof of which we refer to \cite{ab2,agap}. Let $E$ be
a Hilbert space and $f_0$, $G\in C^2(E,\R)$. Let us
denote by $D^2f_0(u)\in{\mathcal L}(E,E')$ the second
Fr\'echet derivative of $f_0$ at $u$. Through the
Riesz Representation Theorem, we can identify
$D^2f_0(u)$ with $f_0''(u)\in {\mathcal L}(E,E)$ given
by $f_0''(u)v={\mathcal K}(D^2 f_0(u)v)$ where
${\mathcal K}:E'\to E$, $\big({\mathcal
  K}(\varphi),\psi\big)_E={}_{\scriptscriptstyle
  E'}\langle\varphi,\psi\rangle_{\scriptscriptstyle E}$, for
  any $\varphi\in E', \psi\in E$. Suppose that $f_0$ satisfies
\begin{align*}
{\rm (a)}& \quad f_0\ \text{\rm has a finite dimensional manifold of critical
  points}\ Z; \\
{\rm (b)}& \quad \text{\rm for all $z\in Z$,  $f_0''(z)$ is a Fredholm operator
  of index $0$;}\\
{\rm (c)}& \quad \text{\rm for all $z\in Z$,  there results $T_zZ=\ker f_0''(z)$.}
\end{align*}
Condition (c) is in fact a nondegeneracy condition which is needed to
apply the Implicit Function Theorem. The inclusion
$T_zZ\subseteq\ker f_0''(z)$ always holds due to the criticality of~$Z$,
so that to prove (c) is enough to prove that $\ker f_0''(z)\subseteq
T_zZ$. \par
Consider the perturbed functional $f_{\e}(u)=f_0(u)-\e G(u)$, and
denote by $\Gamma$ the functional $G\big|_{Z}$. Due to assumptions
(a), (b), and (c), it is possible to prove (see Lemma \ref{l:w}) that
there exists, for $|\e|$ small, a smooth function $w_{\e}(z):\
Z\to(T_zZ)^{\perp}$ such that any critical point $\bar z\in Z$ of the
functional
\[
\Phi_{\e}:\ Z\longrightarrow\R,\qquad \Phi_{\e}(z)=f_{\e}(z+w_{\e}(z))
\]
gives rise to a critical point $u_{\e}=\bar z+w_{\e}(\bar z)$ of
$f_{\e}$; in other words, the perturbed manifold $Z_{\e}=\{z+w_{\e}(z):\ z\in
Z\}$ is a {\em natural constraint} for $f_{\e}$. Moreover $\Phi_{\e}$ admits an expansion of the type
\begin{equation}\label{eq:exp}
\Phi_{\e}(z)=b-\e\Gamma(z)+o(\e)\quad\text{as}\ \e\to 0
\end{equation}
where $b=f_0(z)$ for any $z\in Z$.
\begin{thm}\label{t:abstract}
Let $f_0$ satisfy {\rm (a)}, {\rm (b)}, and {\rm (c)} and assume that $\Gamma$ has a
proper local maximum or minimum point $\bar z$. Then for $|\e|$ small enough,
the functional $f_{\e}$ has a critical point $u_{\e}$ such that
$u_{\e}\to \bar z$ as $\e\to 0$.
\end{thm}
\begin{rem}\label{r:abs}
If $Z_0=\{z\in Z:\ \Gamma(z)=\min_Z\Gamma\}$ is compact, it is still
possible to prove that $f_{\e}$ has a critical point near $Z_0$. The
set $Z_0$ can also consist of local minimum points; the same holds for
maximum points.
\end{rem}
\subsection{The unperturbed problem}

In order to apply the abstract result stated above, we have to prove
that the unperturbed functional
\[
f_0(u)=\frac12\int_{\He^n}|\n_{\He^n}u|^2-\frac1{Q^\star}\int_{\He^n}u_+^{Q^\star},\quad
u\in \SC,
\]
satisfies (a), (b) and (c). Condition (a) clearly
holds; indeed, as remarked above, $f_0$ has a two
dimensional manifold of critical points $Z$, see
(\ref{eq:critman}). Moreover, it is quite standard to
prove that (b) holds. Indeed $f_0\in C^2(\SC,\R)$,
\[
\big(f_0''(u)v,h\big)=(v,h)-(Q^\star-1)\int_{\He^n}u_+^{Q^\star-2}v\,h\quad\forall\,
u,v,h\in\SC,
\]
and for any $z\in Z$ the operator $f_0''(z):\SC\to\SC$
is of the type $I-C_z$, where $I$ is the identity and
$C_z$ is a compact operator and hence $f_0''(z)$ is a
Fredholm operator of index $0$.

Let us now prove (c). The tangent space $T_{z_{\mu,s}}Z$ is given by
\begin{equation*}
T_{z_{\mu,s}}Z=\bigg\{\alpha\frac{\partial z_{\l,t}}{\partial
  \l}\ristretto_{\substack{\l=\mu\\t=s}}+\beta\frac{\partial z_{\l,t}}{\partial t}\ristretto_{\substack{\l=\mu\\t=s}}:\ \alpha,\beta\in\R\bigg\}.
\end{equation*}
\begin{lem}\label{l:nondegeneracy}
For any $\mu>0$, $s\in\R$, there holds
\[
T_{z_{\mu,s}}Z=\ker f_0''( z_{\mu,s}).
\]
\end{lem}
\begin{pf}
As remarked in the previous subsection, it is enough to prove the
inclusion $\ker f_0''(z_{\mu,s})\subseteq T_{z_{\mu,s}}Z$. We can
assume $\mu=1$ and $s=0$, since $\ker f_0''(z_{\mu,s})$ is isomorphic
to $\ker f_0''(z_{1,0})$, due to the invariance of the problem under
dilations and translations along the $t$-axis. If $u\in\ker
f_0''(z_{1,0})$, we have that $u$ is a solution in the weak sense of
$\SC$ of the linearized problem
\[
-\D_{\He^n}u=(Q^\star-1)z_{1,0}^{Q^\star-2}u\quad\text{in}\ \SC.
\]
Due to Lemma \ref{wcs}, $u$ solves the above equation
also in $S^1_0(\He^n)$. It has been proved by
Malchiodi and one of the authors \cite{mu} that any
solution $u$ in $S^1_0(\He^n)$ of such equation is of
the type
\[
u(\xi')=\alpha \frac{\partial \omega_{\l,\xi}}{\partial
  \l}\ristretto_{\substack{\l=1\\\xi=0}}(\xi')+\sum_{i=1}^{2n+1}\nu_i\frac{\partial
  \omega_{\l,\xi}}{\partial
  \xi_i}\ristretto_{\substack{\l=1\\\xi=0}}(\xi'),\quad\quad \xi'\in \He^n,
\]
for some coefficients $\alpha\in\R$ and
$\nu=(\nu_1,\dots,\nu_{2n+1})\in\R^{2n+1}$, where, for
$\l\in\R$ and $\xi\in\He^n$,
$\omega_{\l,\xi}=\l^{-\frac{Q-2}2}\omega\circ\delta_{\l^{-1}}\circ\tau_{\xi^{-1}}$,
namely if $\xi'=(z',t')=(x',y',t')$,
$\xi=(z,t)=(x,y,t)$, then
\begin{align}\label{eq:deco}
u(z',t')=&u(x',y',t')=\alpha\, \frac{\partial \omega_{\l,x,y,t}}{\partial
 \l}\ristretto_{(1,0,0,0)}(x',y',t')+\beta\,\frac{\partial
  \omega_{\l,x,y,t}}{\partial
  t}\ristretto_{(1,0,0,0)}(x',y',t')\notag\\
&+\sum_{i=1}^{n}\gamma_i\,\frac{\partial
  \omega_{\l,x,y,t}}{\partial
  x_i}\ristretto_{(1,0,0,0)}(x',y',t')+\sum_{i=1}^{n}\tau_i\,\frac{\partial
  \omega_{\l,x,y,t}}{\partial
  y_i}\ristretto_{(1,0,0,0)}(x',y',t')
\end{align}
for some $\alpha,\beta\in\R$,
$\gamma=(\gamma_1,\dots,\gamma_{n})\in\R^{n}$, and
$\tau=(\tau_1,\dots,\tau_{n})\in\R^{n}$. We claim that the cylindrical symmetry of $u$ implies that, for any
$i=1,2,\dots,n$, $\gamma_i=\tau_i=0$. We have that
\begin{align*}
\frac{\partial
  \omega_{\l,\xi}}{\partial\l}\ristretto_{\substack{\l=1\\\xi=0}}(z',t')=&-\frac{Q-2}2\,\omega(z',t')+(Q-2)^{{Q}/{2}}\big({t'}^2+(1+|z'|^2)^2\big)^{-\frac{Q+2}4}(1+|z'|^2)|z'|^2\\[-10pt]
&+(Q-2)^{{Q}/{2}}\big({t'}^2+(1+|z'|^2)^2\big)^{-\frac{Q+2}4}{t'}^2,\\[5pt]
\frac{\partial
  \omega_{\l,\xi}}{\partial\xi}\ristretto_{\substack{\l=1\\\xi=0}}(z',t')
=&\frac{(Q-2)^{{Q}/{2}}}{\big({t'}^2+(1+|z'|^2)^2\big)^{\frac{Q+2}4}}\big((1+|z'|^2)x'-t'y',(1+|z'|^2)y'+t'x',t'/2\big).
\end{align*}
Therefore $\frac{\partial
  \omega_{\l,\xi}}{\partial\l}\lower3pt\hbox{\big|}_{(1,0)}$ and
$\frac{\partial \omega_{\l,x,y,t}}{\partial
  t}\lower3pt\hbox{\big|}_{(1,0,0,0)}(x',y',t')=\frac12(Q-2)^{{Q}/{2}}\big({t'}^2+(1+|z'|^2)^2\big)^{-\frac{Q+2}4}t'
$
are cylindrically
  symmetric functions. If $u$ is cylindrical, in view of the
  cylindrical symmetry of $\frac{\partial
  \omega_{\l,x,y,t}}{\partial \l}$ and $\frac{\partial
  \omega_{\l,x,y,t}}{\partial t}$, from
(\ref{eq:deco}) we deduce that
\[
\sum_{i=1}^{n}\gamma_i\,\frac{\partial
  \omega_{\l,x,y,t}}{\partial
  x_i}\ristretto_{(1,0,0,0)}(z',t')+\sum_{i=1}^{n}\tau_i\,\frac{\partial
  \omega_{\l,x,y,t}}{\partial
  y_i}\ristretto_{(1,0,0,0)}(z',t')
\]
must be cylindrical, hence $
h(x',y',t')=\sum_{i=1}^{n}\gamma_i\,\big((1+|z'|^2)x_i'-t_i'y_i'\big)+\tau_i
\,\big((1+|z'|^2)y_i'+t_i'x_i'\big) $ must be
cylindrical. From
\[
h(\underbrace{0,\dots,{\mathop{1}\limits^i}\dots,0}_{x'},\underbrace{0,\dots,0}_{y'},0)=h(\underbrace{0,\dots,-{\mathop{1}\limits^i}\dots,0}_{x'},\underbrace{0,\dots,0}_{y'},0)
\]
it follows that $2\gamma_i=-2\gamma_i$ and hence $\gamma_i=0$ for any
$i=1,\dots,n$. In the same way $\tau_i=0$ for any $i=1,\dots,n$. The
claim is thereby proved. As a consequence, we have that $\ker
f_0''(z_{1,0})$ is contained in $T_{z_{1,0}}Z$.
\end{pf}

\subsection{Study of $\Gamma$ and proof of Theorems \ref{t:pert2} and \ref{t:pert3}}

In our case, the reduced functional $\Gamma$ is given by
\begin{align*}
\Gamma(\mu,s)&=\frac1
{Q^\star}\int_{\He^n}k(z,t)z^{Q^{\star}}_{\mu,s}(z,t)\,dz\,dt=\frac1{Q^\star}\int_{\He^n}k(\mu z,\mu^2 t+s)\omega^{Q^{\star}}(z,t)\,dz\,dt\\
&=\frac{\gamma_n}{Q^\star}\int\limits_{\substack{0<r<\infty\\t\in\R}}k(\mu r,\mu^2
t+s)\omega^{Q^{\star}}(r,t)r^{2n-1}\,dr\,dt,
\end{align*}
where $\gamma_n$ is the measure of the unit
$(2n-1)$-sphere. The function $\Gamma$ can be extended with continuity
to $\mu=0$ by setting
\begin{equation}\label{eq:gamma1}
\Gamma(0,s)=\frac{\gamma_n}{Q^\star}k(0,s)\int\limits_{\substack{0<r<\infty\\t\in\R}}\omega^{Q^{\star}}(r,t)r^{2n-1}\,dr\,dt=b_0\,k(0,s)
\end{equation}
where $b_0=\frac1{Q^\star}\int_{\He^n}\omega^{Q^\star}$. By the
Dominated Convergence Theorem it follows that
\begin{equation}\label{eq:gammainfty}
\Gamma(\infty):=\lim_{(\mu,s)\to\infty}\Gamma(\mu,s)=b_0k(\infty).
\end{equation}
Moreover, if $k=\bar k\circ F^{-1}$ with $\bar k$ a smooth function on
the sphere ${\mathbb S}^{2n+1}$, we have that \textsf{}
\begin{align}\label{eq:gamma2}
&D_{\mu}\Gamma(\mu,s)=\frac1{Q^\star}\int_{\He^n}\big[\n_zk(\mu z,\mu^2
t+s)\cdot z+2\partial_t k(\mu z,\mu^2 t+s)\mu
t\big]\omega^{Q^\star}(z,t)\,dz\,dt\notag\\
&D_{\mu}\Gamma(0,s)=\frac1{Q^\star}\int_{\He^n}\big[\n_zk(0,s)\cdot
z\big]\omega^{Q^\star}(z,t)\,dz\,dt=0
\end{align}
and
\begin{align}\label{eq:gamma3}
D^2_{\mu,\mu}\Gamma(0,s)&=\frac1{Q^\star}\int_{\He^n}\bigg[\sum_{i,j=1}^n\frac{\partial^2k}{\partial
  x_i\partial x_j}(0,s)x_i x_j+\sum_{i,j=1}^n\frac{\partial^2k}{\partial
  y_i\partial y_j}(0,s)y_i y_j+2\partial_t
k(0,s)t\bigg]\notag\\
&\qquad\qquad\qquad\cdot\omega^{Q^\star}(x,y,t)\,dx\,dy\,dt\notag\\
&=\frac1{2n Q^\star}\Delta_{x,y}
  k(0,s)\int_{\He^n}|z|^2\omega^{Q^\star}(z,t)\,dz\,dt.
\end{align}

\begin{pfn}{Theorem \ref{t:pert2}}
Assumption (\ref{eq:5}) and (\ref{eq:gammainfty})
imply that $\Gamma(\infty)<\Gamma(\mu,s)$ whereas
assumption (\ref{eq:5}) and (\ref{eq:gamma1}) imply
that $\sup_{\sigma}\Gamma(0,\sigma)<\Gamma(\mu,s)$
hence $\Gamma$ must have a compact set of global
maximum points in the interior of the half-plane
$\{(\mu, s):\ \mu>0\}$. From Remark
\ref{r:abs} we get the conclusion. In the case of
(\ref{eq:oppositesense}) we have $\Gamma(\infty)>\Gamma(\mu,s)$ and
$\inf_{\sigma}\Gamma(0,\sigma)>\Gamma(\mu,s)$ hence $\Gamma$
 must have a compact set of global
minimum points.
\end{pfn}
\begin{pfn}{Theorem \ref{t:pert3}}
The assumptions of Theorem \ref{t:pert3} imply, in
view of \eqref{eq:gamma1}, \eqref{eq:gammainfty},
\eqref{eq:gamma2} and (\ref{eq:gamma3}), that
$\Gamma(0,\bar s)=\max_{\sigma}\Gamma(0,\sigma)\geq
\Gamma(\infty)$ and $D_{\mu}\Gamma(0,\bar s)=0$,
$D^2_{\mu,\mu}\Gamma(0,\bar s)>0$ hence $\Gamma$ must
have a compact set of global maximum points in the
interior of the half-plane $\{(\mu, s):\
\mu>0\}$. The conclusion follows from Remark
\ref{r:abs}.
\end{pfn}

\subsection{Study of $\Phi_{\e}$}

To prove Theorem \ref{t:pert4}, the
study of the functional $\Gamma$ is not sufficient since in this case
$\Gamma$ may be constant even if $k$ is a non-constant function. This fact leads to a loss of information, being the first order expansion
(\ref{eq:exp}) not enough to deduce the existence of critical points
of $\Phi_{\e}$ from the existence of critical points
of $\Gamma$. Therefore we need to study directly the function
$\Phi_{\e}$ which in our case is a function of the two variables $(\mu,s)\in\R^+\times\R$. \par
For $\mu>0$ and
$s\in\R$, let us define the map
\[
{\mathcal U}_{\mu,s}:\ \SC\longrightarrow\SC,\quad{\mathcal
  U}_{\mu,s}(u)(z,t)=\mu^{-\frac{Q-2}2}u\bigg(\frac{z}{\mu},\frac{t-s}{\mu^2}\bigg).
\]
It is easy to check that
$
\big\|{\mathcal U}_{\mu,s}(u)\big\|_{\SC}=\|u\|_{\SC}
$
, for any $u\in\SC$, $\mu>0$, and $s\in\R$,
and that $f_0=f_0\circ{\mathcal U}_{\mu,s}$. Moreover we have that
$
\big({\mathcal U}_{\mu,s}\big)^{-1}={\mathcal
  U}_{\mu^{-1},-\mu^{-2}s}=\big({\mathcal U}_{\mu,s}\big)^{t}
$
where $\big({\mathcal U}_{\mu,s}\big)^{t}$ denotes the adjoint of
${\mathcal U}_{\mu,s}$. Differentiating the identity
$f_0=f_0\circ{\mathcal U}_{\mu,s}$ we observe that
\[
f_0'=\big({\mathcal U}_{\mu,s}\big)^{-1}\circ f_0'\circ {\mathcal
  U}_{\mu,s}
\]
and
\begin{equation}\label{eq:fnd}
f_0''(u)=\big({\mathcal U}_{\mu,s}\big)^{-1}\circ f_0''\big({\mathcal
  U}_{\mu,s}(u)\big)\circ{\mathcal U}_{\mu,s},\quad\forall\, u\in\SC.
\end{equation}
Clearly we have
$
{\mathcal U}_{\mu,s}:\
T_{\omega}Z\longrightarrow T_{z_{\mu,s}}Z\quad\text{and}\quad
{\mathcal U}_{\mu,s}:\
\big(T_{\omega}Z)^{\perp}\longrightarrow\big(T_{z_{\mu,s}}Z\big)^{\perp}.
$
Because of nondegeneracy, the self adjoint Fredholm operator
$f_0''(\omega)$ maps $\SC$ into $ \big(T_{\omega}Z)^{\perp}$ and
$f_0''(\omega)\in {\mathcal
  L}\big((T_{\omega}Z)^{\perp}\big)$. Moreover (\ref{eq:fnd})
implies that
\begin{equation}\label{eq:6}
\|f_0''(\omega)^{-1}\|_{{\mathcal
  L}((T_{\omega}Z)^{\perp})}
=\|f_0''(z)^{-1}\|_{{\mathcal
  L}((T_{z}Z)^{\perp})}\quad\forall\, z\in Z.
\end{equation}
\begin{lem}\label{l:w}
Assume $k\in L^{\infty}(\He^n)$. Then there exist constants $\e_0,C>0$ and a smooth function
\[
w=w(\mu,s,\e):\quad (0,+\infty)\times \R \times (-\e_0,\e_0)\ \longrightarrow\
\SC
\]
such that for any $\mu>0$, $s\in\R$, and $\e\in(-\e_0,\e_0)$
\begin{align}
\label{eq:w1}
w(\mu,s,\e)\ \text{ is orthogonal to }\ T_{z_{\mu,s}}Z\\
\label{eq:w2}
f_{\e}'\big(z_{\mu,s}+w(\mu,s,\e)\big)\in T_{z_{\mu,s}}Z\\
\label{eq:w3}
\|w(\mu,s,\e)\|\leq C\,|\e|.
\end{align}
\end{lem}
\begin{pf}
Since it will be useful in the sequel, we write the complete proof of the lemma
which follows the proofs of analogous results of \cite{ab2,am}. Let us define
\begin{align*}
H:\ &(0,+\infty)\times\R\times \SC\times\R\times\R\times\R\longrightarrow\SC\times\R\times\R\\
&(\mu,s,w,\alpha_1,\alpha_2,\e)\longmapsto \big(f_\e'(z_{\mu,s}+w)-\alpha_1
\dot{\xi}_{\mu,s}-\alpha_2
  \dot{\zeta}_{\mu,s},(w,\dot{\xi}_{\mu,s}),(w,\dot{\zeta}_{\mu,s})\big),
\end{align*}
where $\dot{\xi}_{\mu,s}$ (resp. $\dot{\zeta}_{\mu,s}$) denotes
the normalized tangent vector
$\frac{\partial}{\partial\mu} z_{\mu,s}$
(resp. $\frac{\partial}{\partial s}z_{\mu,s}$). If
$H(\mu,s,w,\alpha_1,\alpha_2,\e)=0$ then $w$ satisfies (\ref{eq:w1})-(\ref{eq:w2}) and
$H(\mu,s,w,\alpha_1,\alpha_2,\e)=0$ if and only if
$(w,\alpha_1,\alpha_2)$ is a fixed point for the map $F_{\mu,s,\e}$
defined as
\[
F_{\mu,s,\e}(w,\alpha_1,\alpha_2):= -\bigg(\frac{\partial H}{\partial(w,\alpha_1,\alpha_2)}(\mu,s,0,0,0,0)\bigg)^{-1}
H(\mu,s,w,\alpha_1,\alpha_2,\e) + (w,\alpha_1,\alpha_2).
\]
To prove the existence of $w$ satisfying (\ref{eq:w1}) and
(\ref{eq:w2}) it is enough to prove that $F_{\mu,s,\e}$ is a contraction in
some ball $B_\rho(0)$, with $\rho=\rho(\e)>0$
independent of $z \in Z$, whereas the regularity of $w(\mu,s,\e)$
follows from the Implicit Function Theorem. We have that
\[
\bigg(\frac{\partial H}{\partial(w,\alpha_1,\alpha_2)}(\mu,s,0,0,0,0)\bigg)(w,\beta_1, {\beta}_2)=
\big(f_0''(z_{\mu,s})w - \beta_1 \dot{\xi}_{\mu,s}-  {\beta}_2
\dot{\zeta}_{\mu,s},(w,\dot{\xi}_{\mu,s}),(w,\dot{\zeta}_{\mu,s}) \big).
\]
>From (b) we deduce that $\big(\frac{\partial H}{\partial(w,\alpha_1,\alpha_2)}
(\mu,s,0,0,0)\big)$ is an injective Fredholm operator of index zero,
hence it is invertible and
\begin{align*}
&\bigg(\frac{\partial
  H}{\partial(w,\alpha_1,\alpha_2)}(\mu,s,0,0,0,0)\bigg)^{-1}(w,\beta_1, {\beta}_2)\\
&\quad=
\Big(\beta_1 \dot{\xi}_{\mu,s}+  {\beta}_2
\dot{\zeta}_{\mu,s}+f_0''(z_{\mu,s})^{-1}\big(w-(w,\dot{\xi}_{\mu,s})\dot{\xi}_{\mu,s}-(w,\dot{\zeta}_{\mu,s})\dot{\zeta}_{\mu,s}\big),-(w,\dot{\xi}_{\mu,s}),-(w,\dot{\zeta}_{\mu,s})\Big ).
\end{align*}
In view of (\ref{eq:6}), we have that
$
\big\|\big(\frac{\partial H}{\partial(w,\alpha_1,\alpha_2)}(\mu,s,0,0,0,0)\big)^{-1}\big\| \leq
\max\big(1,\|(f_0''(z_{\mu,s}))^{-1}\|\big)
= \max\big(1,\|(f_0''(\omega))^{-1}\|\big) .
$
Set $C_*=\max\big(1,\|(f_0''(\omega))^{-1}\|\big)$. For any $(w,\alpha_1,\alpha_2)\in B_\rho(0)$ we have that
\begin{align}\label{eq:fmse}
&\|F_{\mu,s,\e}(w,\alpha_1,\alpha_2)\| \leq C_* \|f_\e'(z_{\mu,s}+w)-f_0''(z_{\mu,s})w\|\notag \\
&\qquad\qquad\leq C_*\int_0^1 \|f_0''(z_{\mu,s}+tw)-f_0''(z_{\mu,s})\|\,dt\cdot \|w\| +
C_*|\e| \|G'(z_{\mu,s}+w)\|\notag \\
&\qquad\qquad\leq C_*\int_0^1 \|f_0''(\omega+t\,{\mathcal U}_{\mu,s}^{-1}(w))-f_0''(\omega)\|\,dt\cdot\|w\| +
C_*|\e| \|G'(z_{\mu,s}+w)\|\notag\\
&\qquad\qquad\leq C_*\rho \sup_{\|w\|\le \rho }\|f_0''(\omega+w)-f_0''(\omega)\|
+ C_*|\e| \sup_{\|w\|\le \rho }\|G'(z_{\mu,s}+w)\| .
\end{align}
For $(w_1,\alpha_1,{\beta}_1),(w_2,\alpha_2, {\beta}_2) \in B_\rho(0)$
\begin{align}\label{eq:fmse1}
&\frac{\|F_{\mu,s,\e}(w_1,\alpha_1,{\beta}_1)-F_{\mu,s,\e}(w_2,\alpha_2, {\beta}_2)\|}{C_*
 \|w_1-w_2\|}\notag\\
&\qquad\qquad\qquad\qquad \leq
\frac{\|f_\e'(z_{\mu,s}+w_1) - f_\e'(z_{\mu,s}+w_2)-f_0''(z_{\mu,s})(w_1-w_2)\|}{\|w_1-w_2\|}\notag \\
&\qquad\qquad\qquad\qquad\leq \int_0^1 \|f_0''(z_{\mu,s}+w_2+t(w_1-w_2))-f_0''(z_{\mu,s})\|
\, dt\notag \\
&\qquad\qquad\qquad\qquad\quad +|\e| \int_0^1 \|G''(z_{\mu,s}+w_2+t(w_1-w_2))\|
\,dt \notag\\
&\qquad\qquad\qquad\qquad\le \sup_{\|w\|\leq 3\rho }\|f_0''(\omega+w)-f_0''(\omega)\| + |\e|
\sup_{\|w\|\le 3 \rho}\|G''(z_{\mu,s}+w)\|.
\end{align}
Choose $\rho_0>$ such that
$
C_* \sup_{\|w\|\le 3\rho_0
  }\|f_0''(\omega+w)-f_0''(\omega)\|<{1}/{2}
$
and $\e_0>0$ such that
\begin{align*}
2\e_0<\Big(\sup_{z \in Z, \|w\|\le 3 \rho_0}\|G''(z+w)\|\Big)^{-1}C_*^{-1} \text{ and }
3\e_0< \Big(\sup_{z \in Z, \|w\|\le \rho_0}\|G'(z+w)\|\Big)^{-1}C_*^{-1}\rho_0.
\end{align*}
With these choices, for any $z_{\mu,s} \in Z$ and
$|\e|<\e_0$ the map $F_{\mu,s,\e}$ maps $B_{\rho_0}(0)$ into itself and is a contraction
there such that
$\|F_{\mu,s,\e}(w_1,\alpha_1,{\beta}_1)-F_{\mu,s,\e}(w_2,\alpha_2, {\beta}_2)\|\leq\l\|w_1-w_2\|$,
where the constant $\l\in(0,1)$ does not depend on $\mu,s,\e$. Therefore
$F_{\mu,s,\e}$ has a unique fixed point $(w(\mu,s,\e),\alpha_1(\mu,s,\e),{\alpha}_2(\mu,s,\e))$ in $B_{\rho_0}(0)$. From (\ref{eq:fmse}) we also infer that $F_{\mu,s,\e}$ maps $B_\rho(0)$ into $B_\rho(0)$, whenever
$\rho\le\rho_0$ and
\[
\rho> 2|\e| \big(\sup_{\|w\|\le \rho}\|G'(z_{\mu,s}+w)\|\big)C_*.
\]
Consequently for the uniqueness of the fixed point we have
\begin{align*}
\|(w(\mu,s,\e),\alpha_1(\mu,s,\e),{\alpha}_2(\mu,s,\e))\| \le  3 |\e| \big(\sup_{\|w\|\le \rho_0}\|G'(z_{\mu,s}+w)\|\big)C_*,
\end{align*}
which gives (\ref{eq:w3}).
\end{pf}

\noindent
We are now interested in the behavior of the function
\[
\Phi_{\e}(\mu,s)=f_{\e}\big(z_{\mu,s}+w(\mu,s,\e)\big)
\]
the critical points of which on $\R^+\times\R$ give rise to critical
points of $f_{\e}$ on $\SC$, as remarked in Subsection
\ref{sec:abs}. In particular we will prove the following proposition.
\begin{pro}\label{p:limitphi}
Assume that $k$ is cylindrically symmetric and continuous on
$\ov{\He^n}$. Then for any $\bar s\in\R$ there holds
\[
\begin{array}{rll}
\text{\rm(i)}&\lim_{(\mu,s)\to(0,\bar s)}\Phi_{\e}(\mu,s)&=f_0(\omega)\big(1+\e
k(0,\bar s)\big)^{-\frac{Q-2}2}\\
\text{\rm(ii)}&\lim_{(\mu,s)\to \infty}\Phi_{\e}(\mu,s)&=f_0(\omega)\big(1+\e
k(\infty)\big)^{-\frac{Q-2}2}.
\end{array}
\]
\end{pro}
\noindent
\begin{rem}\label{r:continuity}
Thanks to the above proposition we can extend $\Phi_{\e}$ to the axis
$\{\mu=0\}$ by setting
\[
\Phi_{\e}(0,s):= f_0(\omega)\big(1+\e
k(0, s)\big)^{-\frac{Q-2}2}
\]
and to infinity by setting
\[
\Phi_{\e}(\infty):=f_0(\omega)\big(1+\e
k(\infty)\big)^{-\frac{Q-2}2}
\]
thus obtaining a continuous function on the compactified half-plane
 $\{(\mu,s): \ \mu\geq0\}\cup\{\infty\}$.
\end{rem}
\noindent For $\mu>0,\ s\in\R$, let us consider the functional
$f_{\e}^{\mu,s}=f_{\e}\circ {\mathcal U}_{\mu,s}$ i. e.
\begin{equation*}
f_{\e}^{\mu,s}(u)=\frac12 \int_{\He^n}|\n_{\He^n}u|^2\,dz\,dt-\frac1{Q^\star}
\int_{\He^n}\big(1+\e k(\mu z,\mu^2 t+s)\big)u_+^{Q^\star}\,dz\,dt.
\end{equation*}
There results $\big(f_{\e}^{\mu,s}\big)'=\big({\mathcal
  U}_{\mu,s}\big)^{-1}\circ f_{\e}'\circ {\mathcal U}_{\mu,s}$ and $\big(f_{\e}^{\mu,s}\big)''(u)=\big({\mathcal U}_{\mu,s}\big)^{-1}\circ
f_{\e}''(u)\circ {\mathcal U}_{\mu,s}$, for any $u\in\SC$.
Let us consider the map $\wtilde H^{\mu,s}:\
\SC\times\R\times\R\times\R\longrightarrow\SC\times\R\times\R$, $(w,\alpha_1,\alpha_2,\e)\longmapsto \big((f^{\mu,s}_\e)'(\omega+w)-\alpha_1
\dot{\xi}_0-\alpha_2
  \dot{\zeta}_0,(w,\dot{\xi}_{0}),(w,\dot{\zeta}_{0})\big)$
where $\dot{\xi}_{0}$ (resp. $\dot{\zeta}_{0}$) is
normalized tangent vector $\frac{\partial}{\partial\mu}
z_{\mu,s}\lower3pt\hbox{\big|}_{\mu=1,s=0}$
\Big(resp. $\frac{\partial}{\partial s}
z_{\mu,s\lower3pt\hbox{\big|}_{\mu=1,s=0}}$\Big). We have that
\begin{align*}
\frac{\partial \wtilde H^{\mu,s}}{\partial(w,\alpha_1,\alpha_2)}(0,0,0,0)=\frac{\partial  H}{\partial(w,\alpha_1,\alpha_2)}(\mu,s,0,0,0,0)\ristretto_{\substack{\mu=1\\s=0}}
\end{align*}
hence $\frac{\partial \wtilde
  H^{\mu,s}}{\partial(w,\alpha_1,\alpha_2)}(0,0,0,0)$ is invertible
  and
$
\big\|\big(\frac{\partial \wtilde
  H^{\mu,s}}{\partial(w,\alpha_1,\alpha_2)}(0,0,0,0)\big)^{-1}\big\| \leq C_*.
$
The map
\[
F^{\mu,s}_{\e}(w,\alpha_1,\alpha_2):= -\bigg(\frac{\partial \wtilde H^{\mu,s}}{\partial(w,\alpha_1,\alpha_2)}(0,0,0,0)\bigg)^{-1}
\wtilde H^{\mu,s}(w,\alpha_1,\alpha_2,\e) + (w,\alpha_1,\alpha_2)
\]
satisfies
\begin{align*}
\|F^{\mu,s}_{\e}(w,\alpha_1,\alpha_2)\| \leq C_* \|f_\e'(z_{\mu,s}+{\mathcal
  U}_{\mu,s}(w))-f_0''(z_{\mu,s}){\mathcal U}_{\mu,s}(w)\|
\end{align*}
and
\begin{align*}
&\frac{\|F^{\mu,s}_{\e}(w_1,\alpha_1,{\beta}_1)-F^{\mu,s}_{\e}(w_2,\alpha_2, {\beta}_2)\|}{C_*
 \|w_1-w_2\|}\\
&\qquad\qquad \leq \frac{\|f_\e'(z_{\mu,s}+{\mathcal U}_{\mu,s}(w_1))-f_\e'(z_{\mu,s}+{\mathcal U}_{\mu,s}(w_2)))-f_0''(z_{\mu,s})({\mathcal U}_{\mu,s}(w_1-w_2))\|}{\|{\mathcal U}_{\mu,s}(w_1)-{\mathcal U}_{\mu,s}(w_2)\|}
\end{align*}
which imply, in view of (\ref{eq:fmse}) and (\ref{eq:fmse1}), that
$F^{\mu,s}_{\e}$ is a contraction in the same ball where
$F_{\mu,s,\e}$ is a contraction (see the proof of Lemma \ref{l:w}).
Hence $F^{\mu,s}_{\e}$ has a fixed
point $(w_{\e}^{\mu,s},\alpha_{1,\e}^{\mu,s},{\alpha}_{2,\e}^{\mu,s})$
such that $\wtilde
H^{\mu,s}(w_{\e}^{\mu,s},\alpha_{1,\e}^{\mu,s},{\alpha}_{2,\e}^{\mu,s},\e)=0$.
>From the uniqueness of the fixed point of $F^{\mu,s}_{\e}$ and from
the fact that
$
(f_{\e}^{\mu,s})'(\omega+({\mathcal U}_{\mu,s})^{-1}w(\mu,s,\e))\in
T_{\omega}Z
$
and
$
({\mathcal U}_{\mu,s})^{-1}w(\mu,s,\e)\in
\big(T_{\omega}Z\big)^{\perp}
$,
it follows that $w^{\mu,s}_{\e}=({\mathcal
  U}_{\mu,s})^{-1}\big(w(\mu,s,\e)\big)$, where $w(\mu,s,\e)$ is given in Lemma
  \ref{l:w}. Assume now that $k$ is continuous on $\ov{\He^n}$ and fix
  $s\in\R$. Let us consider the functional
\[
f^{0,s}_{\e}=\frac12 \int_{\He^n}|\n_{\He^n}u|^2-\frac1{Q^\star}(1+\e
k(0,s)) \int_{\He^n} u_+^{Q^\star}.
\]
For $w_{\e}^{0,s}=(t_{\e}(s)-1)\omega$ where $t_{\e}(s)=(1+\e
k(0,s))^{-\frac{Q-2}4}$ we have
$
(f_{\e}^{0,s})'(\omega+w_{\e}^{0,s})=0
$
and
$
(f_{\e}^{0,s})(\omega+w_{\e}^{0,s})=(1+\e
k(0,s))^{-\frac{Q-2}2}\big(\frac 12-\frac
1{Q^\star}\big)\int_{\He^n}\omega ^{Q^\star}
$
and hence $\wtilde H^{0,s}(w^{0,s}_{\e},0,0,\e)=0$ where
\[
\wtilde
H^{0,s}(w,\alpha_1,\alpha_2,\e)=\big((f_{\e}^{0,s})'(\omega+w)-\alpha_1\dot\xi_0-\alpha_2\dot\zeta_0,
(w,\dot\xi_0), (w,\dot\zeta_0)\big).
\]
We have that $\frac{\partial \wtilde
  H^{0,s}}{\partial(w,\alpha_1,\alpha_2)}(0,0,0,0)=\frac{\partial
  H}{\partial(w,\alpha_1,\alpha_2)}(\mu,s,0,0,0,0)\lower3pt\hbox{\big|}_{\substack{\mu=1\\s=0}}$
  and hence $(w_{\e}^{0,s},0,0)$ is a fixed point of the map
\[
F^{0,s}_{\e}(w,\alpha_1,\alpha_2)=-\bigg(\frac{\partial \wtilde
  H^{0,s}}{\partial(w,\alpha_1,\alpha_2)}(0,0,0,0)\bigg)^{-1}\wtilde H^{0,s}(w,\alpha_1,\alpha_2,\e)+(w,\alpha_1,\alpha_2).
\]
It is easy to check that $F^{0,s}_{\e}$ is a contraction in some ball
of radius $O(|\e|)$. Hence $w^{0,s}_{\e}$ is the unique
fixed point of $F^{0,s}_{\e}$ in such a ball. \\
Let us also set $w^{\infty}_{\e}=(t^{\infty}_{\e}-1)\omega$, $t^{\infty}_{\e}=(1+\e
k(\infty))^{-\frac{Q-2}4}$.
\begin{lem}\label{l:convw}
For any $\bar s\in\R$ there holds
\begin{align}
w^{\mu,s}_{\e}&\to w^{0,\bar s}_{\e},\quad\text{as}\ (\mu,s)\to
(0,\bar s),\label{eq:convW0}\\
w^{\mu,s}_{\e}&\to w^{\infty}_{\e},\quad\text{as}\ (\mu,s)\to
(0,\infty).\label{eq:convWinf}
\end{align}
\end{lem}
\begin{pf}
We have that
\begin{align*}
&\|F_{\e}^{\mu,s}(w_{\e}^{0,\bar s},0,0)-F_{\e}^{0,\bar
  s}(w_{\e}^{0,\bar s},0,0)\|\leq C_* \big\|\wtilde H^{\mu,s}(w_{\e}^{0,\bar s},0,0,\e)-\wtilde H^{0,\bar s}(w_{\e}^{0,\bar s},0,0,\e)\big\| \\
&\quad \leq
C_{*}\|(f_{\e}^{\mu,s})'(\omega+w_{\e}^{0,\bar s})-(f_{\e}^{0,\bar s})'(\omega+w_{\e}^{0,\bar s}
)\|=C_{*}\|(f_{\e}^{\mu,s})'(t_{\e}(\bar s)\omega)-(f_{\e}^{0,\bar s})'(t_{\e}(\bar s)\omega)\|.
\end{align*}
Since by (\ref{eq:sobolev}) and the H\"older inequality
\begin{align*}
&\bigg|\big((f_{\e}^{\mu,s})'(t_{\e}(\bar s)\omega)-(f_{\e}^{0,\bar s})'(t_{\e}(\bar s)\omega),
v\big)\bigg| =\bigg|\int_{\He^n}\e\big[ k(\mu z,\mu^2 t+s)-
k(0,\bar s)\big](t_{\e}(\bar s)\omega)^{Q^\star-1}v\bigg|\\
&\qquad\qquad\qquad \leq
S^{-1/2}\|v\|\bigg[\int_{\He^n}\e^{\frac{Q^\star}{Q^\star-1}}|k(\mu z,\mu^2 t+s)-
k(0,\bar s)|^{\frac{Q^\star}{Q^\star-1}}t_{\e}(\bar s)^{Q^\star}\omega^{Q^\star}\bigg]^{\frac{Q^
\star-1}{Q^\star}}
\end{align*}
we have that, by the Dominated Convergence Theorem,
\begin{equation*}
\|(f_{\e}^{\mu,s})'(t_{\e}(\bar s)\omega)-(f_{\e}^{0,\bar s})'(t_{\e}(\bar s)\omega)\|\leq c
\bigg[\int_{\He^n}|k(\mu z,\mu^2 t+s)-
k(0,\bar s)|^{\frac{Q^\star}{Q^\star-1}}\omega^{Q^\star}\bigg]^{\frac{Q^
\star-1}{Q^\star}}\!\!\! \mathop{\longrightarrow}\limits_{(\mu,s)\to
  (0,\bar s)}0.
\end{equation*}
Therefore
\begin{equation}\label{eq:convF}
F_{\e}^{\mu,s}(w^{0,\bar s}_{\e},0,0)\to
F_{\e}^{0,\bar s}(w^{0,\bar s}_{\e},0,0),\quad\text{as}\ (\mu,s)\to
  (0,\bar s).
\end{equation}
Since $F_{\e}^{\mu,s}$
is a contraction with a contraction factor $0<\l<1$ independent of
$\mu,s$, and $\e$ we have that
\begin{align*}
&\|w^{\mu,s}_{\e}-w^{0,\bar s}_{\e}\|\leq \|(w^{\mu,s}_{\e},\alpha^{\mu,s}_{1,\e},
{\alpha}^{\mu,s}_{2,\e})-(w^{0,\bar s}_{\e},0,0)\|=
\|F^{\mu,s}_{\e}(w^{\mu,s}_{\e},\alpha^{\mu,s}_{1,\e},{\alpha}^{\mu,s}_{2,\e})
-F^{0,\bar s}_{\e}(w^{0,\bar s}_{\e},0,0)\|\\
&\qquad\leq
\|F^{\mu,s}_{\e}(w^{\mu,s}_{\e},\alpha^{\mu,s}_{1,\e},{\alpha}^{\mu,s}_{2,\e})-F^{\mu,s}_{\e}
(w^{0,\bar s}_{\e},0,0)\|+\|F^{\mu,s}_{\e}(w^{0,\bar s}_{\e},0,0)-F^{0,\bar s}_{\e}
(w^{0,\bar s}_{\e},0,0)\|\\
&\qquad\leq \l\|w^{\mu,s}_{\e}-w^{0,\bar s}_{\e}\|+ \|F^{\mu,s}_{\e}(w^{0,\bar s}_{\e},0,0)
-F^{0,\bar s}_{\e}(w^{0,\bar s}_{\e},0,0)\|
\end{align*}
and hence from (\ref{eq:convF}) we obtain
(\ref{eq:convW0}). The proof of (\ref{eq:convWinf}) is analogous.
\end{pf}

\begin{pfn}{Proposition \ref{p:limitphi}}
By definition of $\Phi_{\e}$ and $f^{\mu,s}_{\e}$, we have that
\begin{equation}\label{eq:phi1}
\Phi_{\e}(\mu,s)=f_{\e}\big(z_{\mu,s}+w(\mu,s,\e)\big)=f^{\mu,s}_{\e}(\omega+w_{\e}^{\mu,s}).
\end{equation}
Moreover
\begin{align*}
f^{\mu,s}_{\e}(\omega+w_{\e}^{\mu,s})&=f^{\mu,s}_{\e}(\omega+w_{\e}^{\mu,s})-f^{0,\bar s}_{\e}(\omega+w_{\e}^{\mu,s})+f^{0,\bar s}_{\e}(\omega+w_{\e}^{\mu,s})\\
& =\frac{\e}{Q^\star}\int_{\He^n}\big(k(0,\bar s)-k(\mu z,\mu^2 t+s)\big)(\omega+w_{\e}^{\mu,s})_+^{Q^\star}\,dz\,dt+f^{0,\bar s}_{\e}(\omega+w_{\e}^{\mu,s})
\end{align*}
hence by (\ref{eq:convW0}) and the Dominated Convergence Theorem
\begin{equation}\label{eq:phi2}
f^{\mu,s}_{\e}(\omega+w_{\e}^{\mu,s})\
\mathop{\longrightarrow}\limits_{(\mu,s)\to(0,\bar s)}\
f^{0,\bar s}_{\e}(\omega+w_{\e}^{0,\bar s}).
\end{equation}
On the other hand we have that
\begin{align}\label{eq:phi3}
f^{0,\bar s}_{\e}(\omega+w_{\e}^{0,\bar s})=f^{0,\bar s}_{\e}(t_{\e}(\bar s)\omega)=f_0(\omega)\big(1+\e k(0,\bar s)\big)^{-\frac{Q-2}2}.
\end{align}
>From (\ref{eq:phi1}), (\ref{eq:phi2}), and (\ref{eq:phi3}), (i)
follows. In an analogous way, using (\ref{eq:convWinf}), it is easy
to prove (ii).
\end{pfn}

Thanks to Proposition \ref{p:limitphi}
 it is now quite easy to give
some conditions on $k$ in order to have critical points of
$\Phi_{\e}$. In particular the knowledge of $k$ on the axis $\{(0,s):\
s\in\R\}$ and at $\infty$ gives exact informations about the behavior
of $\Phi_{\e}$ on the axis $\{(0,s):\
s\in\R\}$ and at $\infty$.

\begin{pfn}{Theorem \ref{t:pert4}} As remarked in Subsection
  \ref{sec:abs}, it is enough to prove that $\Phi_{\e}(\mu,s):\
  \R^+\times\R\to\R$ has a critical   point. In view Proposition
  \ref{p:limitphi}, $k(0,t)=k(\infty)$ $\forall\,t\in\R$ implies that
  $\Phi_{\e}(0,t)=\Phi_{\e}(\infty)$ $\forall\,t\in\R$. Hence, either
  $\Phi_{\e}$ is constant (and we have infinitely many critical
  points) or it has a global maximum or minimum point $(\bar \mu,\bar
  s)$,~$\bar\mu>0$. In any case, $\Phi_{\e}$ has a critical point which
  provides a solution to $({\mathcal P}_{\e})$.
\end{pfn}

\appendix
\section*{Appendix}
\setcounter{section}{1}
\setcounter{equation}{0}
\setcounter{thm}{0}

In the first part of this appendix, we prove some technical lemmas
about the properties of cylindrically symmetric functions of the Folland-Stein Sobolev space $S^1_0(\He^n)$.
\begin{rem}\label{r:cyl-grad}
 If $u(z,t)=\ttt{u}(|z|,t)$ , $v(z,t)=\ttt{v}(|z|,t)$ are in
 $\SC$, then $\ttt{u},\ttt{v}\in
 H^1_{\loc}(\{(r,t)\in\R^2|r>0\})$.
 Moreover the following formula holds a.e.
\begin{equation}\label{cg}
  \langle\nh u, \nh v\rangle (z,t)=(\de_r\ttt{u} \de_r\ttt{v}
  +4r^2\de_t\ttt{u} \de_t\ttt{v})(|z|,t).
\end{equation}
\end{rem}

\begin{pf}
It is straightforward to verify that formula \eqref{cg} holds
for smooth functions. In order to extend it to general
$u,v\in \SC$, we choose two sequences of cylindrically
symmetric functions $\phi_k,\psi_k\in C_0^{\infty}(\He^n)$,
converging in $S_0^1$ to $u,v$, respectively. By a
cylindrical change of coordinates, it is then easy to see
that $\ttt{\phi}_k, \ttt{\psi}_k$ are Cauchy sequences in
$H^1(\Omega)$ for every $\Omega\subset\subset
]0,\infty[\times\R$. Since moreover $\phi_k\to u$, $\psi_k\to
v$, pointwise a.e. (up to subsequences), the limits of
$\ttt{\phi}_k, \ttt{\psi}_k$ in $H^1(\Omega)$  are
necessarily $\ttt{u}, \ttt{v}$. As a consequence, formula
\eqref{cg}  (which we know to hold for $\phi_k,\psi_k$)
extends to $u,v$ by means of the a.e. pointwise convergences
$\nh\phi_k\to\nh u$, $\nh\psi_k\to\nh v$,
$(\de_r\ttt{\phi}_k, \de_t\ttt{\phi}_k) \to (\de_r\ttt{u},
\de_t\ttt{u})$,  $(\de_r\ttt{\psi}_k, \de_t\ttt{\psi}_k) \to
(\de_r\ttt{v}, \de_t\ttt{v})$.
\end{pf}

\begin{lem}\label{wcs}
  Let $K(z,t)=\ttt{K}(|z|,t)$, with $\ttt{K}$ bounded and
  locally H\"older continuous in $]0,\infty[\times\R$, and
  let $u\in\SC$ be a (nonnegative)
   weak solution of $-\DHH u=Ku^{Q^*-1}$ in
  the $\SC$-sense (i.e. with respect to $\SC$-test functions).
  Then $u$ is a weak solution of the same equation in the
  $S_0^1(\He^n)$-sense. An analogous result also holds for
  weak solutions of the equation $-\DHH u
  =(Q^*-1)U_{\mu,s}^{Q^*-2}u$.
\end{lem}

\begin{pf}
Let us prove the statement related to the equation
$-\DHH u=Ku^{Q^*-1}$ (the same proof works also for
the other equation). Recalling Remark
\ref{r:cyl-grad}, for every test function
$\ttt{\phi}\in C_0^\infty(]0,\infty[\times\R)$ we have
\begin{equation*}
\begin{split}
  c_n\int_{]0,\infty[\times\R}
  \ttt{K}\ttt{u}^{Q^*-1}\ttt{\phi} r^{2n-1}drdt
&
  =\int_{\He^n}
  Ku^{Q^*-1}\phi
  =\int_{\He^n}
  \langle\nh u,\nh \phi\rangle \\
  &
  =c_n\int_{]0,\infty[\times\R}
  (\de_r\ttt{u} \de_r\ttt{\phi}
  +4r^2\de_t\ttt{u} \de_t\ttt{\phi}) r^{2n-1}drdt
\end{split}
\end{equation*}
where $u(z,t)=\tilde u(|z|,t)$. Hence $\ttt{u}$ is a weak solution of
the elliptic equation
$$-\de_r(r^{2n-1}\de_r\ttt{u})-
  \de_t(4r^{2n+1}\de_t\ttt{u})
  = r^{2n-1}\ttt{K}\ttt{u}^{Q^*-1}$$
in $]0,\infty[\times\R$.
    Now, using a classical bootstrap elliptic argument, one
    can easily see that $\ttt{u}\in C^2(]0,\infty[\times\R)$.
As a consequence $u\in C^2(\{(z,t)\in\He^n|z\neq 0\})$
is a classical solution of the equation $-\DHH
u=Ku^{Q^*-1}$ in $\{z\neq 0\}$ and we can argue as
follows. Let us fix a test function $\phi\in
C_0^\infty(\He^n)$ with support contained in
$\R^{2n}\times[-T,T]$ and let us set
$\Omega_\e=\{(z,t)\in\He^n: |z|<\e, |t|<T\}$,
$\delta_\e=\{(z,t)\in\He^n: |z|=\e, |t|\le T\}$. Let
us also choose a vanishing sequence of positive
numbers $\e_k$ such that
$$\int_{\delta_{\e_k}} |\nh u|^2 dH_{2n}=o(\tfrac{1}{\e_k}),
\qquad\text{as }k\to \infty$$ (such a sequence does
exist since $|\nh u|\in L^2(\He^n)$). Then, setting
\[
A=
\begin{pmatrix}
I_{n} & 0     &  2y     \\
0     & I_{n} & -2x     \\
2y    & -2x   & 4|z|^{2}
\end{pmatrix},
\]
we have by the Divergence Theorem
\begin{equation*}
\begin{split}
& \bigg|\int_{\He^n\setminus\Omega_{\e_k}} (\langle\nh
u,\nh \phi\rangle-
  Ku^{Q^*-1}\phi)\bigg|=
  \bigg|\int_{\He^n\setminus\Omega_{\e_k}} {\rm div}(\phi A\n u)\bigg|\\
  &\quad
    =\bigg|\int_{\delta_{\e_k}} \langle \phi A\n u,\n(-z)\rangle dH_{2n}\bigg|
  =\bigg|\int_{\delta_{\e_k}}\phi \langle  \nh u,\nh(z)\rangle
  dH_{2n}\bigg|\\
  &\quad
  \le c
  \int_{\delta_{\e_k}} |\nh u| dH_{2n}
\le c
 \bigg( \int_{\delta_{\e_k}} |\nh u|^2 dH_{2n}\bigg)^{1/2}\e_k^{(2n-1)/2}
 =o(\e_k^{n-1}),
\quad\text{as }k\to \infty.
\end{split}
\end{equation*}
Since $\big|\int_{\He^n\setminus\Omega_{\e_k}}
(\langle\nh u,\nh\phi\rangle-
Ku^{Q^*-1}\phi)\big|\to\big|\int_{\He^n}(\langle\nh u,
\nh \phi\rangle-  Ku^{Q^*-1}\phi)\big|$, this proves
that $\int_{\He^n} \langle\nh u,\nh \phi\rangle
=\int_{\He^n} Ku^{Q^*-1}\phi$ holds for every $\phi\in
C_0^\infty(\He^n)$  and thus for every $\phi\in
S_0^1(\He^n)$.
\end{pf}

\begin{lem}\label{l:A1}
Let $u_m$ be a sequence weakly converging in $\SC$ to
some function $u\in\SC$. Then (up to subsequences)
$u_m\to u$ in $L^q(C)$ for any set $C$ of the type
$\{z: 0<c_1\leq |z|\leq c_2\}\times [-c_3,c_3]$ (and
hence on any compact set away from the axis $\{z=0\}$)
and for any $1\leq q<+\infty$.
\end{lem}
\begin{pf}
Let $C=\{z:\ c_1\leq |z|\leq c_2\}\times [-c_3,c_3]$.
>From Remark \ref{r:cyl-grad} we have that for any
function $w\in\SC$
\begin{equation}\label{eq:A1}
|\n _{\He^n}w|^2=|\partial_r w|^2+4r^2|\partial_t w|^2.
\end{equation}
Let now $u_m$ be a sequence weakly converging to $u$ in $\SC$ (and so
in $S^1_0(\He^n)$ and in $L^{Q^\star}(\He^n)$). Thanks to
(\ref{eq:A1}) we can write
\begin{align*}
{\rm const}&\geq \int_{C}|\n_{\He^n} u_m|^2=\gamma_n\int_{\substack{c_1\leq r\leq
    c_2\\[2pt]|t|\leq c_3}}\big(|\partial_r u_m|^2+4r^2|\partial_t
    u_m|^2\big)r^{2n-1}\,dr\,dt\\
&\geq \gamma_n\min\{1, 4c_1^2\}c_1^{2n-1}\int_{[c_1,c_2]\times[-c_3,c_3]}
\big(|\partial_r u_m|^2+|\partial_t
    u_m|^2\big)\,dr\,dt
\end{align*}
and analogously for the $L^2$ norm. Hence $u_m(r,t)$
is bounded in $H^1\big([c_1,c_2]\times[-c_3,c_3]
\big)$ which is compactly embedded in
$L^q\big([c_1,c_2]\times[-c_3,c_3] \big)$ for any
$1\leq q<+\infty$. Therefore, up to a subsequence,
$u_m(r,t)\to u(r,t)$ in
$L^q\big([c_1,c_2]\times[-c_3,c_3] \big)$.
Consequently, we get that
\begin{align*}
\int_C |u_m-u|^{q}&=\gamma_n\int_{\substack{c_1\leq r\leq
    c_2\\[2pt]|t|\leq c_3}}|u_m(r,t)-u(r,t)|^{q}r^{2n-1}\,
    dr\,dt\\
    &\leq \gamma_n c_2^{2n-1}\int_{[c_1,c_2]\times[-c_3,c_3]}|u_m(r,t)-u(r,t)|^{q}\,
    dr\,dt\longrightarrow 0.
\end{align*}
Lemma \ref{l:A1} is thereby established.
\end{pf}

Let us now state the P. L. Lions
concentration-compactness principle in $\He^n$. Since
the proof does not present further difficulties with
respect to the euclidean case (see \cite{lions1} and
\cite{lions2}), we omit it. Let
$\ov{\He}^n=\He^n\cup\{\infty\}$ be the
compactification of $\He^n$. Let us denote by
${\mathcal M}(\ov{\He}^n)$ the Banach space of finite
signed regular Borel measures on $\ov{\He}^n$, endowed
with the total variation norm. In view of the Riesz
representation theorem, the space ${\mathcal
M}(\ov{\He}^n)$ can be identified with the dual of the
Banach space $C(\ov{\He}^n)$. We say that a sequence
of measures $\mu_m$ weakly converges to $\mu$ in
${\mathcal M}(\ov{\He}^n)$ if for any $f\in
C(\ov{\He}^n)$ (i.e. continuous on $\He^n$ with finite
limit at $\infty$)
\[
\int_{\He^n}f\, d\mu_m\longrightarrow \int_{\He^n}f\, d\mu.
\]
In this case we will use the notation $\mu_m\ \conm\  \mu$.
\begin{thm}\label{t:cc}
{\bf (Concentration-compactness)} Let $\{u_m\}$ be a
sequence weakly converging to $u$ in $S^1_0(\He^n)$.
Then, up to subsequences,
\begin{align*}
(i)\ \  &\text{$|\n_{\He^n} u_m|^2$ weakly converges
in ${\mathcal M}(\ov{\He}^n)$ to a nonnegative measure
  $\mu$,\quad\qquad}\\
(ii)\ \ &\text{$|u_m|^{Q^\star}$ weakly converges in
${\mathcal M}(\ov{\He}^n)$ to a
  nonnegative measure $\nu$.}
\end{align*}
Moreover there exist an at most countable index set
$J$, a sequence $(z_j,t_j)\in\He^n$,
$\nu^j,\nu^{\infty}\in(0,\infty)$ such that
\begin{align*}
&\nu=|u|^{Q^\star}+\sum_{j\in J}\nu^j\delta_{(z^j,t^j)}+\nu^{\infty}\delta_{\infty}\\
&\mu\geq |\n_{\He^n}u|^2+\sum_{j\in J} S
(\nu^j)^{\frac
  2{Q^\star}}\delta_{(z^j,t^j)}+S(\nu^{\infty})^{\frac
  2{Q^\star}}\delta_{\infty}.
\end{align*}
\end{thm}

\end{document}